\documentclass[11pt]{amsart}
\usepackage{amssymb}
\usepackage{amsmath}
\usepackage{comment}
\usepackage{latexsym}
\usepackage{amsfonts}
\usepackage{graphicx}
\usepackage[normalem]{ulem}
\usepackage[utf8]{inputenc}
\usepackage[bookmarks=true,pdfstartview=FitH, pdfborder={0 0 0},
colorlinks=true,citecolor=blue, linkcolor=blue, hypertexnames=false]{hyperref}
\newtheorem{theorem}{Theorem}

\newtheorem{proposition}[theorem]{Proposition}
\newtheorem{remark}[theorem]{Remark}

\protected\def\vts{
  \ifmmode
    \mskip0.5\thinmuskip
  \else
    \ifhmode
      \kern0.08334em
    \fi
  \fi
}

\input{tcilatex}

\begin{document}

\title[A Partial Generalization of Hantzsche's Theorem]{A Partial Generalization of Hantzsche's Theorem and a Correction}

\author[B. Chow]{Bennett Chow$^{\vts\text{a}}$}
\author[M. Freedman]{Michael Freedman$^{\vts\text{b}}$}
\author[H. Shin]{Henry Shin$^{\vts\text{c}}$}
\author[Y. Zhang]{Yongjia Zhang$^{\vts\text{d}}$}
\date{}

\address{\footnotesize{$^{\text{a}\vts}$\emph{Department of Mathematics, University of California, San Diego, California 92093, USA.}} }

\address{\footnotesize{$^{\text{b}\vts}$\emph{Department of Mathematics, Harvard University, CMSA, Cambridge, Massachusetts 02139, USA.}} }

\address{\footnotesize{$^{\text{c}\vts}$\emph{8800 Lombard Pl., \#1411, San Diego, CA 92122, USA.}} }

\address{\footnotesize{$^{\text{d}\vts}$\emph{School of Mathematical Sciences, Shanghai Jiao Tong University, Shanghai 200240, China.}} }

\begin{abstract}
This note corrects an error in the proof of Proposition 13 in \cite{ChowFreedmanShinZhang2020Advances} and simultaneously establishes a more general result.
We prove that if $M $ is a compact 
connected oriented
$4$-manifold
with connected boundary $\partial M$, 
and if an unbounded number of disjoint copies of $M$ embed topologically and locally flatly in the interior of a compact $4$-manifold $N,$ then $\operatorname{Tor}H_1(\partial M;\mathbb{Z})$ is a direct double, i.e., $\operatorname{Tor}H_1(\partial M;\mathbb{Z})\cong A \oplus A$, with the linking pairing vanishing identically on the first summand, i.e., the linking pairing is split metabolic.
This partially generalizes Hantzsche's theorem stating that the linking pairing for a closed $3$-manifold that embeds in $S^4$ is hyperbolic.
\end{abstract}

\maketitle

\section{A partial generalization of Hantzsche's theorem}

In this note we begin with what needs to be corrected and proceed to state the most general theorem that our homological methods permit.\smallskip 

Proposition 13 in \cite{ChowFreedmanShinZhang2020Advances} 
is equivalent to 
the following statement:

\begin{proposition}\label{prop: Direct Double orig prop 13}
Let $M$ be a compact connected oriented $4$-manifold with connected boundary $\partial M$ diffeomorphic to a spherical space form $S^3/\Gamma$.
If there exists a compact $4$-manifold $N$ containing an unbounded number of
pairwise disjoint copies of $M$ in its interior,
then $H_1(S^3/\Gamma;\mathbb{Z} )$ is a direct double, i.e., isomorphic to $A \oplus A$ for some finite abelian group $A$.
\end{proposition}

As pointed out to us by Max Hallgren, the last sentence in the proof of Proposition 13 in \cite{ChowFreedmanShinZhang2020Advances} contains an erroneous inference.
In particular, in \cite{ChowFreedmanShinZhang2020Advances}, we have proved that 
there exists a finite abelian group $A$ for which
$$H_1(S^3/\Gamma;\mathbb{Z}_p )\cong H_1(S^3/\Gamma;\mathbb{Z} )\otimes \mathbb{Z}_p
\cong (A\otimes \mathbb{Z}_p) \oplus (A\otimes \mathbb{Z}_p)
$$ 
is a direct double for all primes $p$,
and that $|H_1(S^3/\Gamma;\mathbb{Z} )| =|A|^2$.
However, one cannot deduce from these facts alone that $H_1(S^3/\Gamma; \mathbb{Z})$ is a direct double. 
In this paper, we complete the proof of Proposition 13 in \cite{ChowFreedmanShinZhang2020Advances}, while at the same time providing Theorem \ref{thm:07}, which appears
to be the limit of the present method.\smallskip 

Recall that for any closed oriented $3$-manifold $X$, the torsion linking pairing
\[
L_X :
\operatorname{Tor} H_1(X;\mathbb{Z})
\times
\operatorname{Tor} H_1(X;\mathbb{Z})
\longrightarrow
\mathbb{Q}/\mathbb{Z}
\]
is defined as follows: for $x,y \in \operatorname{Tor} H_1(X;\mathbb{Z})$,
choose $k>0$ with $k x = 0$, represent $x$ by a $1$-cycle $a$, choose a
$2$-chain $c$ with $\partial c = k a$, and define
\[
L_X (x,y)
=
\frac{1}{k}\,(c \cdot b)
\in \mathbb{Q}/\mathbb{Z},
\]
where $b$ represents $y$ and $c \cdot b$ denotes algebraic intersection number.
This is well-defined and yields a nonsingular symmetric bilinear pairing.

\begin{theorem}\label{thm:07}
Let $M$ be a compact
connected oriented
$4$-manifold with connected boundary
$\partial := \partial M$ (no homology-sphere assumption), and let $N$ be a compact 
$4$-manifold (with or without boundary). 
If, for every $n \in \mathbb{N}$, there exist disjoint locally flat
codimension-zero topological embeddings $\bigsqcup^{n} M \hookrightarrow \operatorname{int}(N)$, then the torsion subgroup 
\[
B \ :=\ \operatorname{Tor} H_1(\partial\vts; \mathbb{Z})
\]
is a direct double. More precisely, if
\[
A \ :=\ \operatorname{Tor} H_1(M; \mathbb{Z}),
\]
then $B \cong A \oplus A$.\smallskip  

Moreover, the torsion linking pairing
$L$ on $B \cong A \oplus A$ 
vanishes on the first summand $A \oplus 0$, which we may identify with
$\ker\bigl(\operatorname{Tor} H_1(\partial\vts;\mathbb{Z}) \to H_1(M;\mathbb{Z})\bigr)$, i.e., $L$ is split metabolic. 
\end{theorem}

\begin{remark}
    Theorem \ref{thm:07} partially generalizes Hantzsche's theorem \cite{Hantzsche}, which states that
if a closed $3$-manifold $X$ embeds in $S^4,$ then $\operatorname{Tor}(H_1(X;\mathbb{Z}) )  \cong A\oplus A$ for some finite abelian group $A$, and the linking pairing of $X$ is hyperbolic, vanishing on both $A$-summands.\smallskip 

Indeed, any closed $3$-manifold $X$ that embeds in $S^4$ separates $S^4$ into two components. Let $M$ be the closure of one of the components, so that $\partial M=X.$
Theorem \ref{thm:07} says that if an unbounded number of disjoint copies of $M$ embed as pairwise disjoint locally flat codimension-zero topological submanifolds of the interior of some compact $4$-manifold $N$, much of Hantzsche's conclusion still holds.
\end{remark}

\begin{proof}
Set
\[
K := H_2(M,\partial\vts;\mathbb{Z}),\qquad
\widehat B := H_1(\partial\vts;\mathbb{Z}),\qquad
C := H_1(M;\mathbb{Z}),
\]
\[
A := \operatorname{Tor} H_1(M;\mathbb{Z})=\operatorname{Tor} C,
\qquad \text{and} \qquad 
B := \operatorname{Tor} H_1(\partial\vts;\mathbb{Z})=\operatorname{Tor}\widehat B.
\]

We first claim that for every coefficient ring
\[
R=\mathbb{Z}\quad\text{or}\quad R=\mathbb{Z}/m\mathbb{Z}\ \ (m>0),
\]
the map
\[
i_*:H_k(\partial\vts;R)\to H_k(M;R)
\]
is surjective for $k=1,2$.\smallskip 

Indeed, for an embedding $e_n:\bigsqcup^n M\hookrightarrow N$, let
\[
M^n:=e_n\bigl(\bigsqcup^n M\bigr),\qquad
C_n:=N\setminus \operatorname{int}(M^n).
\]
Then $\partial^n:=\partial M^n\cong \bigsqcup^n \partial$.
Since the embeddings are locally flat codimension-zero embeddings, $\partial M^n$
is bicollared in $N$. Hence there exist open neighborhoods $U$ of $M^n$ and
$V$ of $C_n$ such that $U$ deformation retracts onto $M^n$, $V$ deformation
retracts onto $C_n$, and $U\cap V$ deformation retracts onto $\partial^n$.
Applying Mayer--Vietoris to the open cover $N=U\cup V$ gives
\[
H_k(\partial^n;R) \xrightarrow{(i_*^{\oplus n},j_*)}
H_k(M^n;R)\oplus H_k(C_n;R)\longrightarrow H_k(N;R)
\qquad (k= 1,2).
\]
Projecting to the first factor shows that
\[
\operatorname{coker}(i_*^{\oplus n})
\]
is a quotient of $\operatorname{coker}(i_*^{\oplus n},j_*)$, hence a subquotient of $H_k(N;R)$.\smallskip 

If $\operatorname{coker}(i_*)\neq 0$, then
\[
\operatorname{coker}(i_*^{\oplus n})
\cong \operatorname{coker}(i_*)^{\oplus n}.
\]
For $R=\mathbb{Z}$, the minimal number of generators of this group grows unboundedly with $n$, contradicting the fact that it is a subquotient of the fixed finitely generated group $H_k(N;\mathbb{Z})$.
For $R=\mathbb{Z}/m\mathbb{Z}$, the group $\operatorname{coker}(i_*)$ is a nonzero finite $R$-module, so
\[
\bigl|\operatorname{coker}(i_*^{\oplus n})\bigr|
=
\bigl|\operatorname{coker}(i_*)\bigr|^n
\]
grows unboundedly with $n$, contradicting the fact that it is a subquotient of the fixed finite group $H_k(N;\mathbb{Z}/m\mathbb{Z})$.
Thus $i_*$ is surjective for $k=1,2$ and for all such coefficient rings $R$.\smallskip 

Since $M$ and $\partial$ are connected, the map
\[
H_0(\partial\vts;\mathbb{Z}) \longrightarrow H_0(M;\mathbb{Z})
\]
is an isomorphism. Together with the surjectivity of
\[
i_*:H_1(\partial\vts;\mathbb{Z}) \to H_1(M;\mathbb{Z}),
\]
the long exact sequence of the pair $(M,\partial)$ implies that
\[
H_1(M,\partial\vts;\mathbb{Z})=0.
\]
Moreover, since
\[
i_*:H_2(\partial\vts;\mathbb{Z}) \to H_2(M;\mathbb{Z})
\]
is surjective, the map
\[
H_2(M;\mathbb{Z}) \to H_2(M,\partial\vts;\mathbb{Z})
\]
is zero. Hence the long exact sequence yields the short exact sequence
\begin{equation}\label{eq:full-ses-general}
0 \to K \to \widehat B \to C \to 0.
\end{equation}

Next, fix $m>0$ and apply the same argument with coefficients $\mathbb{Z}/m\mathbb{Z}$.
Again using the long exact sequence of the pair, we obtain
\[
0 \to H_2(M,\partial\vts;\mathbb{Z}/m\mathbb{Z})
\to H_1(\partial\vts;\mathbb{Z}/m\mathbb{Z})
\to H_1(M;\mathbb{Z}/m\mathbb{Z})
\to 0.
\]
 Thus the universal coefficient theorem gives
\[
H_2(M,\partial\vts;\mathbb{Z}/m\mathbb{Z})
\cong H_2(M,\partial\vts;\mathbb{Z})\otimes \mathbb{Z}/m\mathbb{Z}
\cong K/mK.
\]
Also,
\[
H_1(\partial\vts;\mathbb{Z}/m\mathbb{Z})\cong \widehat B/m\widehat B,
\qquad
H_1(M;\mathbb{Z}/m\mathbb{Z})\cong C/mC.
\]
Therefore, for every $m>0$, we have a short exact sequence
\begin{equation}\label{eq:modm-ses}
0 \to K/mK \to \widehat B/m\widehat B \to C/mC \to 0.
\end{equation}

Recall that a subgroup $S$ of an abelian group $G$ is called \emph{pure} if for every
$n \in \mathbb{Z}_{>0}$ and every $a \in S$, if there exists $x \in G$ such
that $nx = a$, then there exists $y \in S$ such that $ny = a$.\smallskip 

Exactness on the left in \eqref{eq:modm-ses} is equivalent to
\[
K\cap m\widehat B = mK
\qquad\text{for every }m>0.
\]
Thus $K$ is a pure subgroup of the finitely generated abelian group $\widehat B$.
A standard structure-theorem argument then implies that $K$ is a direct summand of $\widehat B$.
Hence \eqref{eq:full-ses-general} splits:
\[
\widehat B \cong K \oplus C.
\]

Taking torsion gives
\[
B=\operatorname{Tor}\widehat B
\cong \operatorname{Tor}K \oplus \operatorname{Tor}C.
\]
By Lefschetz duality,
\[
K=H_2(M,\partial\vts;\mathbb{Z})\cong H^2(M;\mathbb{Z}),
\]
and the cohomological universal coefficient theorem yields
\[
0\to \operatorname{Ext}^1(H_1(M;\mathbb{Z}),\mathbb{Z})
\to H^2(M;\mathbb{Z})
\to \operatorname{Hom}(H_2(M;\mathbb{Z}),\mathbb{Z})
\to 0.
\]
Since $\operatorname{Hom}(H_2(M;\mathbb{Z}),\mathbb{Z})$ is torsion-free, it follows that
\[
\operatorname{Tor}K
\cong \operatorname{Ext}^1(H_1(M;\mathbb{Z}),\mathbb{Z})
\cong \operatorname{Tor}H_1(M;\mathbb{Z})
= A.
\]
Therefore
\[
B \cong A \oplus A.
\]

Finally, from \eqref{eq:full-ses-general} we have
\[
\operatorname{Tor}K
=
\ker\bigl(\operatorname{Tor}H_1(\partial\vts;\mathbb{Z})
\to H_1(M;\mathbb{Z})\bigr).
\]
Let $i:\partial M \hookrightarrow M$ denote the inclusion. 
Under Lefschetz duality for \(M\) and Poincar\'e duality for \(\partial M\), the
map
\[
K=H_2(M,\partial M;\mathbb{Z})\longrightarrow H_1(\partial M;\mathbb{Z})
\]
corresponds, up to the usual sign convention, to the restriction map
\[
i^*:H^2(M;\mathbb{Z})\longrightarrow H^2(\partial M;\mathbb{Z}).
\]
Hence \(\operatorname{Tor}K\) identifies with the image of
\[
i^*:\operatorname{Tor}H^2(M;\mathbb{Z})\longrightarrow
\operatorname{Tor}H^2(\partial M;\mathbb{Z}).
\]

We use the following standard cohomological description of the torsion
linking pairing. If \(X\) is a closed oriented \(3\)-manifold and
\[
\beta_X:H^1(X;\mathbb{Q}/\mathbb{Z})\to H^2(X;\mathbb{Z})
\]
is the Bockstein homomorphism associated to
\[
0\to \mathbb{Z}\to \mathbb{Q}\to \mathbb{Q}/\mathbb{Z}\to 0,
\]
then \(\beta_X\) surjects onto \(\operatorname{Tor}H^2(X;\mathbb{Z})\), and for
\(x,y\in \operatorname{Tor}H^2(X;\mathbb{Z})\) one has, up to the usual overall
sign convention,
\[
L_X(x,y)=\bigl\langle x\smile \widetilde y,[X]\bigr\rangle,
\]
where \(\widetilde y\in H^1(X;\mathbb{Q}/\mathbb{Z})\) satisfies
\(\beta_X(\widetilde y)=y\). This is well defined.

To prove that the subgroup \(\operatorname{Tor}K\) is isotropic, let
\[
x=i^*\alpha,\qquad y=i^*\gamma,
\qquad
\alpha,\gamma\in \operatorname{Tor}H^2(M;\mathbb{Z}).
\]
Choose \(\widetilde\gamma\in H^1(M;\mathbb{Q}/\mathbb{Z})\) such that
\[
\beta_M(\widetilde\gamma)=\gamma.
\]
By naturality of the Bockstein,
\[
\beta_{\partial M}(i^*\widetilde\gamma)=i^*\gamma=y.
\]
Therefore, using the cohomological formula for the torsion linking pairing on
\(\partial M\), we get
\[
\begin{aligned}
L_{\partial M}(x,y)
&=\bigl\langle i^*\alpha\smile i^*\widetilde\gamma,[\partial M]\bigr\rangle\\
&=\bigl\langle i^*(\alpha\smile\widetilde\gamma),[\partial M]\bigr\rangle\\
&=\bigl\langle \alpha\smile\widetilde\gamma,\, i_*[\partial M]\bigr\rangle\\
&=0,
\end{aligned}
\]
because \(i_*[\partial M]=0\in H_3(M;\mathbb{Z})\), equivalently
\([\partial M]=\partial[M,\partial M]\). Thus \(\operatorname{Tor}K\) is
isotropic. Under the splitting
\[
B\cong \operatorname{Tor}K\oplus \operatorname{Tor}C\cong A\oplus A,
\]
it is precisely the first summand. Hence the torsion linking pairing is split
metabolic.
\end{proof}

\begin{remark}
We, of course, wondered if embedding an unbounded number of copies of $M$ in a compact $N$ implied that $M$ actually embeds in a patch, $\mathbb{R}^4$ in either the topological or smooth category. 
Our homological methods offer no hope of resolving this question. 
Nevertheless, we wondered if the full hyperbolicity conclusion of Hantzsche could be extracted for $\partial M$. 
Curiously, the example below shows that the asymptotic size bounds discussed here do not suffice to force hyperbolicity.\smallskip  

For the \(n\)-fold embedding above, let
\[
A^n:=\ker\bigl(\operatorname{Tor}H_1(\partial^n;\mathbb Z)\to H_1(M^n;\mathbb Z)\bigr),
\]
so \(A^n\cong A^{\oplus n}\).
Theorem~\ref{thm:07} identifies one isotropic copy of $A$ (coming from the kernel into $M$). 
The other isotropic summand should be deduced from asymptotic information about the kernel into the complements $C_n$. 
Since $N$ has only finitely generated homology, it is not difficult to argue that this kernel induces isotropic subgroups $S_n \leq  \operatorname{Tor} H_1(\partial^{\,n}; \mathbb{Z})$ with both $|A^n|/|S_n|$ and $|A^n \cap S_n|$ bounded above independently of $n$. 

Again, one wonders if that information is sufficient to conclude that for large $n$
\[
|S_n| = |A^n| \qquad \text{and} \qquad A^n \cap  S_n = \{e\}  ,
\]
in which case the orthogonal sums $L^{\oplus n}$ would necessarily be hyperbolic. 
One might hope that this asymptotic information forces $L^{\oplus n}$ to be hyperbolic for large $n$.
The following example shows that even this is too optimistic.
However, as always with quadratic forms, characteristic $2$ presents an intriguing special case. 
The following example shows that these asymptotic size bounds alone do not
suffice to strengthen ``split metabolic'' to hyperbolic.
\end{remark}

\section{Example of a linking pairing $L$ asymptotically close to hyperbolic}

Let
\[
A = \mathbb{Z}/2\mathbb{Z},
\qquad
T = A \oplus A \cong (\mathbb{Z}/2\mathbb{Z})^2,
\]
and let $L : T \times T \to \mathbb{Q}/\mathbb{Z}$ be the nonsingular symmetric linking pairing given, in the standard basis, by
\[
\Gamma =
\begin{pmatrix}
0 & \frac{1}{2} \smallskip  \\
\frac{1}{2} & \frac{1}{2}
\end{pmatrix}.
\]
Equivalently,
\[
L\big((a, x), (a', x')\big) = \frac{ax' + a'x + xx'}{2} \in \mathbb{Q}/\mathbb{Z},
\qquad
a, a', x, x' \in \mathbb{F}_2.
\]
Then the first summand $A \oplus 0$ is a Lagrangian, since $L\big((1, 0), (1, 0)\big) = 0$. However, $L$ is not hyperbolic:
the two order-$2$ subgroups complementary to $A\oplus 0$ are $\langle(0, 1)\rangle$ and $\langle(1, 1)\rangle$, and
\[
L\big((0, 1), (0, 1)\big) = \frac{1}{2},
\qquad
L\big((1, 1), (1, 1)\big) = \frac{1}{2},
\]
so neither complement is isotropic.\smallskip 

For the $n$-fold orthogonal sum $L^{\oplus n}$, write elements of $T^n$ as $(u, v)$ with $u, v \in (\mathbb{F}_2)^n$. Then
\[
L_n \big( (u, v), (u', v') \big)  = \frac{u \cdot v' + u' \cdot v + v \cdot v'}{2}.
\]
Let
\[
H_n = \left\{ v \in (\mathbb{F}_2)^n : \sum_{i=1}^n v_i = 0 \right\},
\qquad
|H_n| = 2^{n-1},
\]
and define $f_n : H_n \to (\mathbb{F}_2)^n$ by
\[
(f_n(v))_i = \sum_{j<i} v_j.
\]
A direct calculation shows that for all $v, w \in H_n$,
\[
f_n(v) \cdot w + v \cdot f_n(w) = v \cdot w,
\]
hence
\[
G_n := \{(f_n(v), v) : v \in H_n\} \leq  T^n
\]
is isotropic. 
Let
\[
A^n:=A^{\oplus n}\oplus 0 \le T^n.
\]
Since $G_n$ is a graph, we have $G_n \cap A^n=\{0\}$, and
\[
|G_n| = |H_n| = 2^{n-1} = \frac{|A^n|}{2}.
\]
Therefore $L$ is asymptotically hyperbolic in the multiplicative sense: for every fixed $k > 2$,
\[
|G_n| = \frac{|A^n|}{2} > \frac{|A^n|}{k}
\]
for all $n \ge 2$, even though $L$ itself is not hyperbolic.\smallskip

In fact, $L^{\oplus n}$ is not hyperbolic for any $n$.
Indeed, $T^n$ has exponent $2$, so any hyperbolic linking pairing on $T^n$
is alternating. But if $e_1=(1,0,\dots,0)\in (\mathbb{F}_2)^n$, then
\[
L_n\big((0,e_1),(0,e_1)\big)=\frac12,
\]
so $L_n$ is not alternating.\smallskip 

This linking pairing is realized (up to the usual overall sign convention for surgery/linking pairings) by
any closed $3$-manifold obtained by surgery on a $2$-component framed link with linking number $2$ and
self-linkings $2$ and $4$; equivalently, the surgery matrix is
\[
Q =
\begin{pmatrix}
2 & 2 \\
2 & 4
\end{pmatrix},
\qquad
Q^{-1} \equiv
\begin{pmatrix}
0 & \frac{1}{2} \smallskip  \\
\frac{1}{2} & \frac{1}{2}
\end{pmatrix}
\pmod{1}.
\]
\medskip

\textbf{Acknowledgment.} We would like to thank Max Hallgren for pointing out the error in the proof of Proposition 13 in \cite{ChowFreedmanShinZhang2020Advances}.
We would like to thank the referee for helpful suggestions and corrections.
The authors used ChatGPT as a tool for interactive proofreading, error detection, and improving the clarity of several arguments.


\begin{thebibliography}{9}

\bibitem{ChowFreedmanShinZhang2020Advances}Chow, Bennett; Freedman, Michael; Shin, Henry; Zhang,
              Yongjia. \emph{Curvature growth of some $4$-dimensional gradient {R}icci
              soliton singularity models.} Adv. Math.
\textbf{372} (2020), 107303, 17 pp.

\bibitem{FreedmanKrushkal}
Freedman, Michael; Krushkal, Vyacheslav,
\emph{A triple torsion linking form and $3$-manifolds in $S^4$.} arXiv:2506.11941.

\bibitem{Hantzsche}Hantzsche, W. \emph{Einlagerung von {M}annigfaltigkeiten in euklidische {R}\"aume.} Math. Z.
\textbf{43} (1938), 38--58.

\bibitem{Hat}Hatcher, Allen. \emph{Algebraic topology.} Cambridge University Press, Cambridge, 2002. 



\end{thebibliography}
\end{document}